\documentclass{amsart}
\usepackage{amsthm,amsmath,amssymb,graphicx,upref,cite}
\usepackage[line,ps,color,poly,all,xdvi,dvips]{xy}

\CompileMatrices

\newcommand{\R}{\mathbb R}
\newcommand{\Z}{\mathbb Z}
\newcommand{\Q}{\mathbb Q}
\newcommand{\C}{\mathbb C}
\newcommand{\OO}{\mathcal{O}}

\newcommand{\group}[1]{\mathbf{#1}}
\DeclareMathOperator{\SU}{SU}
\DeclareMathOperator{\SL}{SL}
\newcommand{\conj}[1]{\lsp{#1}P_0}
\newcommand{\G}{\Gamma} 

\newcommand{\lsp}[1]{{}^{#1}\!}
\renewcommand{\Im}{\operatorname{Im}}
\newcommand{\QQ}{{\mathcal Q}}

\newcommand{\GG}{{\mathcal G}}

\newcommand{\Par} {{\mathcal P}}
\newcommand{\I} {{\mathcal I}}

\theoremstyle{plain}
\newtheorem{thm}{Theorem}[section]

\newtheorem{prop}[thm]{Proposition}

\theoremstyle{definition}
\newtheorem{defn}[thm]{Definition}

\newtheorem*{rem}{Remark}

\begin{document}


\title{The elliptic points of the Picard modular group over the Gaussian integers}
\author{Dan Yasaki}
\address{Department of Mathematics and Statistics\\Lederle Graduate Research Tower\\ University of Massachusetts\\Amherst, MA 01003-9305}
\email{yasaki@math.umass.edu}
\date{}
\thanks{The original manuscript was prepared with the \AmS-\LaTeX\ macro
system and the \Xy-pic\ package.}
\keywords{Picard modular group, fixed point, elliptic element}
\subjclass[2000]{Primary 11F57; Secondary 53C35}
\begin{abstract}
We explicitly compute the elliptic points and isotropy groups for the action of the Picard modular group over the Gaussian integers on $2$-dimensional complex hyperbolic space.  
\end{abstract}
\maketitle

\bibliographystyle{../../amsplain_initials}

\begin{section}{Introduction}\label{sec:introduction}
Let $D=\group{G}(\R)/K$ be a symmetric space of non-compact type, where $\group{G}$ is a semisimple algebraic group defined over $\Q$.  An arithmetic group $\G\subset \group{G}(\Z)$ acts on $D$ by left translation, and one can study the elliptic points of this action, the points in the interior of $D$ with non-trivial stabilizer.  

One application of this computation is to the study of arithmetic quotients $\Gamma \backslash D$.  The quotient is not smooth in general.  It has orbifold singularities arising from the elliptic elements of $\Gamma$.  An explicit knowledge of the fixed points in $D$ with associated stabilizer groups in $\Gamma$ allow one to study the types of singularities that occur.  

When the $\Q$-rank of $\group{G}$ is 1, one can use a family of exhaustion functions to find elliptic points.  In \cite{Yasrank1}, we define one such family of exhaustion functions.  These exhaustion functions come out of Saper's work on tilings in \cite{Sa}.  In fact, our exhaustion functions are nothing more than the composition of his \emph{normalized parameters} (in the $\Q$-rank 1 case) with the rational root.  In \cite{Yaspicard}, we use the exhaustion functions to construct an explicit deformation retraction of $D$ onto a spine $D_0$ in the case where $\group{G}=\SU(2,1;\Z[i])$ is the Picard modular group over the Gaussian integers. 

In this paper, we use the spine from \cite{Yaspicard} to study the elliptic elements of $\SU(2,1;\Z[i])$.  Section~\ref{sec:general} recalls a general decomposition of $D$ in the $\Q$-rank 1 case.  Section~\ref{sec:background} specializes to the case where $\group{G}=\SU(2,1;\Z[i])$, and Section~\ref{sec:data} looks at the exhaustion functions and spine in this case.  We examine the elliptic points of the complement of the spine in Section~\ref{sec:complement} and the elliptic points in the spine in Section~\ref{sec:spine}.  The elliptic points in $D \setminus D_0$ can be understood in terms of the intersection of $\Gamma$ with a rational parabolic subgroup, and the elliptic points of $D_0$ can be understood in terms of the stabilizers of the cells of $D_0$.

I would like to thank Paul Gunnells for reading an early version of this paper and many helpful comments. I would also like to thank Les Saper for patiently explaining tilings to me.
\end{section}

\begin{section}{General $\Q$-rank 1}\label{sec:general}
In this section we briefly describe a $\Gamma$-invariant decomposition of $D$ into codimension $0$ sets using exhaustion functions.  This construction is described for the general $\Q$-rank 1 case in \cite{Yasrank1}.

\begin{subsection}{Notation}
Let $G=\group{G}(\R)$ be the group of real points of a $\Q$-rank 1 semisimple algebraic group defined over $\Q$.  Let $\Par$ denote the set of proper rational parabolic subgroups of $\group{G}$.  To ease the notation, when there is no risk of confusion, we will use the same roman letter to denote an algebraic group and its group of real points.  
\end{subsection}

\begin{subsection}{The exhaustion functions}
There exists an exhaustion functions $f_P$ for every rational parabolic subgroup $P \subseteq G$.  Since the rational parabolic subgroups correspond to cusps, these functions can be though of as height functions with respect to the various cusps.

The family of exhaustion functions defined above is $\Gamma$-invariant in the sense that 
\begin{equation}\label{eq:invariantfamily}
f_{\lsp{\gamma}P}(z)=f_P(\gamma^{-1}\cdot z)\quad \text{for $\gamma \in \Gamma$.}
\end{equation}
\end{subsection}

\begin{subsection}{Induced decomposition of $D$}
 These exhaustion functions are used to define a decomposition of $D$ into sets $D(\I)$ for $\I \subset \Par$.  For a parabolic subgroup $P$, define $D(P)\subset D$ to be the set of $z \in D$ such that $ f_P(z) \geq  f_Q(z)$ for every $Q \in \Par \setminus \{P\}$.  In other words, $D(P)$ consists of the points that are higher with respect to $P$ than any other cusp.  This gives a decomposition of the symmetric space parameterized by rational parabolic subgroups,
\begin{equation}\label{eq:decomp1}
D=\bigcup_{P \in \Par} D(P).
\end{equation}

More generally, for a subset  $\I \subseteq \Par$, 
\begin{align}
D(\I)&=\bigcap_{P\in \I} D(P) \label{eq:decomp2}\\ 
D'(\I)&=D(\I) \setminus \bigcup_{\tilde{\I} \supsetneq \I} D(\tilde{\I}).
\end{align}
It follows that $D'(\I)\subseteq D(\I)\text{ and } D(\I)=\coprod_{\tilde{\I} \supseteq \I}D'(\tilde{\I})$. 

\begin{defn}
A subset $\I \subset \Par$ is called {\it admissible} if $D(\I)$ is non-empty and {\it strongly admissible} if $D'(\I)$ is non-empty.
\end{defn}

\begin{prop}[{\cite[Proposition~3.7]{Yasrank1}}]\label{prop:tiling}
Let $\mathcal S$ denote the collection of strongly admissible subsets of $\Par$.  Then the symmetric space has a $\Gamma$-invariant decomposition 
\[D = \coprod_{\I \in {\mathcal S}} D'(\I),\]
such that $\gamma \cdot D'(\I)=D'(\lsp{\gamma} \I)$ for all $\gamma \in \Gamma$ and $\I\in \mathcal{S}$.
\end{prop}
\begin{defn}\label{defn:spine}
Given a family of $\Gamma$-invariant exhaustion functions, define a subset $D_0 \subset D$ by
\begin{equation*}
D_0=\coprod_{\substack{\I \in {\mathcal S}\\ |\I| >1}} D'(\I).
\end{equation*} 
\end{defn}
  
\begin{rem}
In fact, there is a $(k-1)$-parameter family of different retractions, where $k$ is the number of $\Gamma$-conjugacy classes of parabolic $\Q$-subgroups.  However, the $\Gamma$-invariance of the exhaustion functions ensure that this choice is immaterial for the purposes of computing elliptic points.  In particular, suppose a point $z\in D$ has stabilizer $\Gamma_z \subset \Gamma$ that is not contained in $\Gamma \cap Q$ for any rational parabolic subgroup $Q$.  Pick a family of exhaustion functions.  Then $z \in D(P)$ for some rational parabolic subgroup $P$.  Since $\Gamma_z$ is not a subgroups of a rational parabolic subgroup, there exists a $\gamma \in \Gamma_z$ such that $\lsp{\gamma}P \neq P$.  Then \eqref{eq:invariantfamily} implies that $f_{\lsp{\gamma}P}(z) = f_P(z)$.  It follows that $z\in D_0$.  
\end{rem}
\end{subsection}
\end{section}

\begin{section}{Background}\label{sec:background}
\begin{subsection}{The unitary group}
Let $G$ be the identity component of the real points of the algebraic group $\group{G}=\SU(2,1)$, realized explicitly as
\[G=\group{G}(\R)=\SU(2,1;\C)=\{g \in \SL(3,\C)\; |\; g^* C g 
= C\},\]
where $C=\begin{bmatrix} 0&0&i\\0&-1&0\\-i&0&0 \end{bmatrix}$.  Alternatively, let $\mathcal{Q}$ be the (2,1)-quadratic form on $\C^3$ defined by 
\begin{equation*}
\mathcal{Q}(u,v)=u^* C v.
\end{equation*}  
Then $G$ is the group of determinant 1 complex linear transformations of $\C^3$ that preserve $\mathcal Q$.  Let $\OO = \Z[i]$ and $\G$ be the arithmetic subgroup $\G = \group{G}(\Z) = G \cap \SL_3(\OO)$. and let $K$ be the maximal compact subgroup $K=G \cap \SU(3)$.

Let $\theta$ denote the Cartan involution given by inverse conjugate transpose and let $K$ be the fixed points under $\theta$.  Then $K$ is the maximal compact subgroup $K=G \cap \SU(3)$.

Because these elements of $\Gamma$ will be used frequently, set once and for all\[w=\begin{bmatrix} 0 & 0 & -1 \\ 0 & 1 & 0 \\ 1 & 0 & 0 \end{bmatrix},
\sigma=\begin{bmatrix} 1 & 1+i & i \\ 0 & 1 & 1+i \\ 0 & 0 & 1 \end{bmatrix}, \check{\sigma}=\begin{bmatrix} 1 & -1+i & i \\ 0 & 1 & 1-i \\ 0 & 0 & 1\end{bmatrix},\] \[\tau= \begin{bmatrix} 1 & 0 & 1 \\ 0 & 1 & 0 \\ 0 & 0 & 1 \end{bmatrix}, \epsilon=\begin{bmatrix} i & 0 & 0 \\ 0 & -1 & 0 \\ 0 & 0 & i \end{bmatrix},\text{ and }\xi=\begin{bmatrix} 1 & -1-i & i \\ 1-i & -1 & 0 \\ 1-i & -1-i & i \end{bmatrix}.\] 
\end{subsection}

\begin{subsection}{The symmetric space}
Let $D=G/K$ be the associated Riemannian symmetric space of non-compact type.  Then $D$ is $2$-dimensional complex hyperbolic space or the complex $2$-ball with the Bergmann metric.  We will put coordinates on $D$ using Langlands decomposition.

Let $P_0 \subset G$ be the rational parabolic subgroup of upper triangular matrices, and fix subgroups $N_0$, $A_0$, and $M_0$:
\begin{align*}
P_0&=\left\{\left.\begin{bmatrix} y\zeta & \beta \zeta^{-2} & \zeta\left(r + i|\beta|^2/2\right)/y \\ 0 & \zeta^{-2} & i{\overline \beta}\zeta/y \\ 0 & 0 & \zeta/y \end{bmatrix} \; \right| \;  \zeta,\beta \in \C,\ |\zeta|=1,\ r \in \R,\ y\in \R_{>0} \right\},\\
N_0&=\left\{\left.\begin{bmatrix} 1 & \beta & r + i|\beta|^2/2 \\ 0 & 1 & i{\overline \beta} \\ 0 & 0 & 1 \end{bmatrix} \; \right| \;  \beta \in \C,\ r \in \R\right\},\\
A_0&=\left\{\left.\begin{bmatrix} y&0&0\\0&1&0\\0&0&1/y \end{bmatrix}\; \right|\;  y \in \R_{>0} \right\},\\
M_0&=\left\{\left.\begin{bmatrix} \zeta&0&0\\0&\zeta^{-2}&0\\0&0&\zeta \end{bmatrix}\; \right|\;  \zeta \in \C, \ |\zeta|=1 \right\}.
\end{align*}
$P_0$ acts transitively on $D$, and every point $z\in D$ can be written as $p \cdot x_0$ for some $p \in P_0$.  Using  Langlands decomposition, there exists $u \in N_0, a \in A_0$, and $m \in M_0$ such that $p=uam$.  Since $M_0 \subset K$, $z$ can be written as $ua \cdot x_0$.  Denote such a point $z=(y,\beta,r)$.  These are also known as \emph{horospherical coordinates}.
\end{subsection}

\begin{subsection}{Parabolic subgroups}
Zink showed that $\Gamma$ has class number 1 \cite{Z}.  Thus $\Gamma \backslash \group{G}(\Q) / \group{P_0}(\Q)$ consists of a single point, and all the rational parabolic subgroups of $G$ are $\Gamma$-conjugate. The rational parabolic subgroups of $G$ are parametrized by the maximal isotropic subspaces of $\C^3$ which they stabilize.  These are $1$-dimensional, and so to each $P \in \Par$, there is an associated primitive, isotropic vector $v_P \in \Z[i]^3$.  (A vector $(n,p,q)^t \in \Z[i]^3$ is \emph{primitive} if $(n,p,q)$ generate $\Z[i]$ as an ideal.)  Similarly, given a primitive, isotropic vector $v$ in $\Z[i]^3$, there is an associated rational parabolic subgroup $P_v$.  Notice, however, that $v_P$ is only well-defined up to scaling by $\Z[i]^*=\{\pm 1, \pm i\}.$  Thus, the vectors $v$ and $\varepsilon v$ will be treated interchangeably for $\varepsilon \in \Z[i]^*$.  If $P=\lsp{\gamma}Q$ for some $\gamma \in \Gamma$, then $v_{P}=\gamma v_{Q}$.  

Unless explicitly mentioned otherwise, the vector $v_P$ will be written as $v_P=(n,p,q)^t$.  The isotropic condition $ \QQ(v_P,v_P)=0$ implies that 
\begin{equation}\label{eq:isotropic}
|p|^2=2\Im(n \overline{q}).
\end{equation}  
In particular, $q\neq 0$ for $P \neq P_0$.  Furthermore, since there are no isotropic 2-planes in $\C^3$,  
\begin{equation*}
\QQ(v_P,v_Q)\neq 0 \quad \text{for $P \neq Q$.}
\end{equation*} 

\end{subsection}
\end{section}
\begin{section}{$\SU(2,1;\Z[i])$ data}\label{sec:data}
We recall some facts shown in \cite{Yaspicard}.
\begin{subsection}{The functions}
We first define an exhaustion function $f_P$ for every rational parabolic subgroup $P \subseteq G$.  

Let $z=(y,\beta,r)\in D$ and $P$ a rational parabolic subgroup of $G$ with associated isotropic vector $(n,p,q)^t$. Then the exhaustion function $f_P$ can be written as
\begin{align}
 f_{0}(z)&\equiv  f_{P_0}(z)=y \label{eq:fP0}\\
 f_P(z)&=\frac{y}{\left(|n-\beta p+\bigl(i|\beta|^2/2-r\bigr)q|^2+y^2|p-i\overline \beta q|^2+y^4|q|^2\right)^{1/2}}\label{eq:fP}.
\end{align}
\end{subsection}

\begin{subsection}{Induced decomposition of $D$}
The strongly admissible sets are completely known.  They correspond to certain configurations of isotopic vectors in $\OO^3$.  In particular, either
\begin{enumerate}
\item $|\I| \leq 5$ and $|\QQ(v_P,v_Q)|^2 \leq 2$ for every $P,Q\in \Par$ or
\item  $|\I| =8$ and $|\QQ(v_P,v_Q)|^2 \leq 4$ for every $P,Q\in \Par$.
\end{enumerate}
The converse is also true.  Every set of primitive, isotropic vectors in $\OO^3$ satisfying the conditions above correspond to a set of rational parabolic subgroups that is strongly admissible.
\end{subsection}

\begin{subsection}{Siegel strip}
\begin{defn}
Let $S\subset D$ be the subset 
\[S=\left\{(y,\beta,r) \in D | -1/2<r \leq 1/2, \beta \in \diamond\right\}, \quad \text{where}\]
$\diamond$ is the square in the complex plane with vertices  $0, \ (1+i)/2, \ i,$ and $(-1+i)/2$. 
\end{defn}
A computation shows that following proposition.
\begin{prop}\label{prop:siegelstrip}
Every point $z \in D$ is conjugate under $\Gamma_{P_0}$ to a point in $S$. 
\end{prop}
\end{subsection}
\end{section}

\begin{section}{Isotropy groups of points in $D\setminus D_0$}\label{sec:complement}
Using Proposition~\ref{prop:tiling}, we divide $D$ into a codimension $1$ piece $D_0$ and a codimension $0$ piece $D \setminus D_0$ and study the fixed points on each separately.  Note that by Proposition~\ref{prop:siegelstrip}, it suffices to look within the intersection of these set with $S$.  Since the structure of $D\setminus D_0$ is easier, we examine this piece first.  

\begin{subsection}{Reduction to computation on $D(P_0) $}
Notice that $D\setminus D_0=\coprod_{P \in \Par} D'(P)$.  Since all of the rational parabolic subgroups are $\Gamma$-conjugate, Proposition~\ref{prop:tiling} implies that it suffices to look in $D(P_0)$. 

\begin{prop}\label{prop:DP0}
Every point in $D$ is a $\G$-translate of a point in $D(P_0)$.  In particular, every point in $D \setminus D_0$ is a $\G$-translate of a point in $D'(P_0)$.
\end{prop}

\begin{prop}\label{prop:Stab0}
The subgroup of $\G$ that stabilizes $D'(P_0)$ is exactly $\G_{P_0}=\G \cap P_0$.
\end{prop}
\begin{proof}
If $\gamma \in \G$ stabilizes $D'(P_0)$, then $\lsp{\gamma} P_0 = P_0$.  Since parabolic subgroups are self-normalizing, $\gamma \in P_0$.
\end{proof}
\end{subsection}

\begin{subsection}{Fixed points in $D'(P_0) \cap S$}  
We compute explicitly the action of $\Gamma_{P_0}$ in coordinates.  With this, we are able to find 

\begin{prop}\label{prop:fix0}
The fixed points in $D'(P_0)\cap S$ are of the form $(y,\beta,r)$ with $\beta \in \{0,(1+i)/2,i\}$.
\end{prop}
\begin{proof}
Let $p=uam=\begin{bmatrix} y\zeta & \beta \zeta^{-2} & \zeta\left(r + i|\beta|^2/2\right)/y \\ 0 & \zeta^{-2} & i{\overline \beta}\zeta/y \\ 0 & 0 & \zeta/y \end{bmatrix}$ as in Section~\ref{sec:background} and let $z_0=(y_0,\beta_0,r_0)$ be a point in $D$.
Then 
\begin{align*}
p\cdot z_0 &=uam \cdot (y_0,\beta_0,r_0)\\
&=ua \cdot (y_0,\zeta^{3}\beta_0,r_0)\\
&=u \cdot (y y_0,y\zeta^{3}\beta_0,y^2 r_0)\\
&=(y y_0,y\zeta^{3}\beta_0+\beta,y^2 r_0+r-\Im(\beta \bar{\beta}_0\zeta^{-3}y)).
\end{align*}
Now suppose $p$ fixes $z_0$.  Then since $y_0>0$, this implies that 
\begin{equation}
 y=1,\quad 
\beta=\beta_0(1-\zeta^3), \quad \text{and} \quad
r=|\beta_0|^2\Im(\zeta^{-3}). 
\end{equation}
Further suppose that $p \in \G_{P_0}$.  Then $r \in \Z,\  \zeta \in \OO^*$, and $\beta \in \OO$ such that $2\mid |\beta|^2$.
It follows that the isotropy group in $\G_{P_0}$ of a point $(y_0,\beta_0,r_0)$ consists of the intersection $\{I,\gamma,\gamma^2,\gamma^3\} \cap \G$, where 
\begin{equation}\label{eq:gamma}
\gamma=\begin{bmatrix} i & -(1+i)\beta_0  & -(1-i)|\beta_0|^2 \\ 0 & -1 & -(1-i){\overline \beta_0} \\ 0 & 0 & i \end{bmatrix}.
\end{equation}
The points of $S$ for which the intersection is non-trivial have $\beta_0\in \{0,(1+i)/2,i\}$.  The intersection has order four for $\beta_0\in \{0,i\}$ and order two for $\beta_0=(1+i)/2$.\end{proof}

\begin{prop}\label{prop:fix0strict}
Every non-trivial isotropy groups of a point in $D\setminus D_0$ is $\G$-conjugate to exactly one of 
\begin{equation*}
\G_1=\langle\epsilon \rangle\cong \Z/4\Z, \quad \G_2=\langle \xi^2 \rangle\cong \Z/4\Z, \quad \text{or} \quad
\G_3=\langle \sigma \epsilon^2 \rangle\cong \Z/2\Z.
\end{equation*}
Furthermore
\begin{align*}
D^1=D^{\G_1}&=\{(y,0,r)\;|\;y >0, r \in \R\},\\
D^2=D^{\G_2}&=\{(y,i,r)\;|\;y >0, r \in \R\}, \quad \text{and}\\
D^3=D^{\G_3}&=\{(y,(1+i)/2,r)\;|\;y >0, r \in \R\},\\
\end{align*}
\end{prop}
\begin{proof}
Let $z\in D\setminus D_0$ be a fixed point.  By Proposition~\ref{prop:DP0}, $z$ is a $\G$-translate of a point in $D'(P_0)$.  The group $\Gamma_{P_0}$ stabilizes $D'(P_0)$.  Propositions~\ref{prop:siegelstrip}~and ~\ref{prop:Stab0} then imply that $z$ is a $\G$-translate of a point in $D'(P_0) \cap S$.  The result then follows from Proposition~\ref{prop:fix0} and evaluating \eqref{eq:gamma} for $\beta_0\in \{0,i,(1+i)/2\}$: 
\begin{align*}
\G_1&=\left \langle\begin{bmatrix} i & 0 & 0 \\ 0 & -1 & 0 \\ 0 & 0 & i \end{bmatrix}\right \rangle=\langle\epsilon \rangle\cong \Z/4\Z,\\
\G_2&=\left \langle \begin{bmatrix} i & 1-i  & -1+i \\ 0 & -1 & 1+i \\ 0 & 0 & i \end{bmatrix}\right \rangle
=\langle \xi^2 \rangle\cong \Z/4\Z, \quad \text{and}\\
\G_3&=\left \langle \begin{bmatrix} -1 & 1+i  & -i \\ 0 & 1 & -1-i \\ 0 & 0 & -1 \end{bmatrix}\right \rangle=\langle \sigma \epsilon^2 \rangle\cong \Z/2\Z.
\end{align*}
\end{proof}
\end{subsection}
\end{section}
\begin{section}{Isotropy groups of points in $D_0$}\label{sec:spine}
In this section, we find the isotropy groups of points on $D_0$.  A fundamental domain for the action of $\G$ on $D_0$ is given in \cite{Yaspicard}.  In particular, $D_0$ is given the structure of a cell complex such that the stabilizer of a cell fixes the cell pointwise.  We mention that this cell structure is just a subdivision of the decomposition given in Definition~\ref{defn:spine}.  

\begin{subsection}{Reduction to computation on representative cells}
The space $D_0$ is given the structure of a cell-complex such that the stabilizer of a cell fixes the cell pointwise.  The cells of $D_0$ fall into $24$ equivalence classes modulo $\G$ consisting of two $3$-cells, seven $2$-cells, nine $1$-cells, and six $0$-cells.  Representatives of the $3$-cells and their boundary faces are shown in Figure~\ref{fig:3cells}.  

\begin{figure}
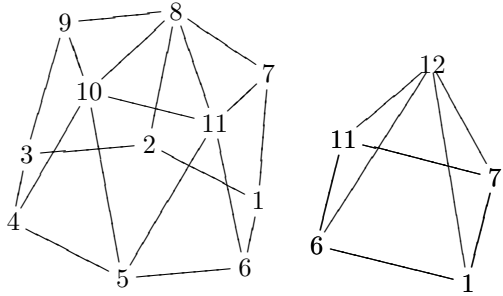

\[\begin{array}{c @{\hspace{0.1in}} c}
{\def\nplus{\ifcase\xypolynode\or 7 \or 8 \or 9 \or 10 \or 11 \fi}
\xy/r.7in/: ="A", +(.2,1.1)="B","A",
{\xypolygon6"C"{~:{(1,-.2):(0,.5)::}~={30}~*{{\xypolynode}}}},
"B",{\xypolygon5"D"{~:{(.8,-.2):(0,.5)::}~*{\nplus}}},
"D2";"D4"**@{-},"D2";"D5"**@{-},
"D5";"C6"**@{-},
"D5";"C5"**@{-},
"D4";"C5"**@{-},
"D4";"C4"**@{-},
"D3";"C3"**@{-},
"D2";"C2"**@{-},
"D1";"C1"**@{-}
\endxy}
&
\def\nplus{\ifcase\xypolynode\or 7 \or 11 \or 6 \or 1 \fi}
\xy/r.7in/: ="A", +(.2,1.1)="B",
"A",
{\xypolygon4{~:{(.8,-.2):(0,.7)::}~*{\nplus}}},
"B"*{12}="C",
"A",
 {\xypolygon4{~:{(.8,-.2):(0,.7)::}~*\nplus~<>{;"C"**@{-}}}}
\endxy
\end{array}\]
\caption{The representative $3$-cells}
\label{fig:3cells}
\end{figure}

\end{subsection}
\begin{subsection}{Fixed points in representative cells}
Since the cell structure is such that the stabilizer of a cell fixes it pointwise, the cells with not-trivial stabilizer form a set of $\G$-representatives of the fixed points that we are looking for.  These are listed in Table~\ref{tab:stab}.
\begin{table}
\caption{Representative cells and their stabilizers}
\label{tab:stab}
\begin{center}
\begin{tabular}{|c|c|c|c|}\hline
Cell & Dimension & Stabilizer & Generators \\\hline
$[1,2,3,4,5,6]$ &2& $\Z/2\Z$ & $\langle \epsilon w \rangle$\\
$[5,11]$ &1& $\Z/2\Z$ & $\langle \sigma \epsilon w \sigma^{-1} \rangle$\\
$[5,6]$&1&$\Z/2\Z$ & $\langle\epsilon w\rangle$\\
$[1,6]$&1&$\Z/2\Z$ &  $\langle \epsilon w \rangle$\\
$[1,2]$&1&$\Z/2\Z$ &  $\langle\epsilon w \rangle$\\
$[2,8]$&1&$\Z/4\Z$ & $\langle \epsilon \rangle$\\
$[1,12]$&1&$\Z/4\Z$ & $\langle \xi^2 \rangle$\\
$[8]$&0&$\Z/12\Z$ & $\langle \tau \epsilon w \rangle$\\
$[2]$&0&$ \Z/2\Z \times \Z/4\Z$ & $\langle \epsilon w , \epsilon \rangle$\\
$[6]$&0&$\Z/2\Z$ & $\langle \epsilon w \rangle $\\
$[1]$&0&$\GG_{31}$\footnotemark & $\langle \epsilon w , \xi^2 \rangle$\\
$[5]$&0&$\mathfrak{S}_3$ &$\langle \epsilon w,\sigma \epsilon^2 \rangle$\\
$[12]$&0&$\Z/8\Z$ &$\langle \xi \rangle$\\\hline
\end{tabular}
\end{center}
\end{table}
\footnotetext{This is the order $32$ group with Hall-Senior number $31$ \cite{Ha} and Magma small group library number $11$.}

Note that $\G_1$, $\G_2$, and $\G_3$ occur as subgroups of the groups in Table~\ref{tab:stab} since $D^i$ contains points with $y$ small and points with $y$ large.  In particular, the surface $D^i$ intersects the spine $D_0$ and so the isotropy group of some point in $D_0$ contains $\Gamma_i$.  

\begin{thm}
The isotropy group of a point in $D$ is $\G$-conjugate to exactly one of the following:
\begin{enumerate}
\item $\G_1=\langle \epsilon \rangle \cong \Z/4\Z$
\item $\G_2=\langle \xi^2 \rangle \cong \Z/4\Z$
\item $\G_3=\langle \sigma \epsilon^2  \rangle \cong \Z/2\Z$
\item $\G_4=\langle \epsilon w\rangle \cong \Z/2\Z$
\item $\G_5=\langle \tau \epsilon w\rangle \cong \Z/12\Z$
\item $\G_6=\langle \epsilon w, \epsilon \rangle \cong \Z/2\Z \times \Z/4\Z $
\item $\G_7=\langle \epsilon w , \xi^2 \rangle \cong \GG_{31}$
\item $\G_8=\langle \epsilon w,\sigma \epsilon^2 \rangle \cong \mathfrak{S}_3$
\item $\G_9=\langle \xi \rangle \cong \Z/8\Z.$ 
\end{enumerate}
\end{thm}
\end{subsection}

\begin{subsection}{Fixed points of $\Gamma_i\ (i \geq 4)$ }
Since the cells are defined using the exhaustion functions, we are able to describe the fixed points explicitly in terms of the exhaustion functions.  

\begin{prop}
Let $D^i=D^{\G_i}$.  Then
\begin{enumerate}
\item $D^4=\left\{(y,\beta,r):y^2+\frac{1}{2}|\beta|^2=1\right\}$\label{it:D4}
\item $D^5=\left\{\left(\sqrt[4]{\frac{3}{4}},0,\frac{1}{2}\right)\right\}$ \label{it:D5}
\item $D^6=\{(1,0,0)\}$
\item $D^7= \left\{\left(\frac{1}{\sqrt{2}},i,0\right)\right\}$
\item $D^8=\left\{\left(\frac{\sqrt{3}}{2},\frac{1+i}{2},0\right)\right\}$
\item $D^9=\left\{\left(\frac{1}{\sqrt[4]{2}},i,\frac{1}{2}\right)\right\}$\label{it:D9}
\end{enumerate}
\end{prop}

\begin{proof}
In particular, $\G_4$ is the stabilizer of $[1,2,3,4,5,6]$.  This is precisely the set $D(\{P_0, \lsp{w}P\})\cap \{r=0\}$.  Item \eqref{it:D4} then follows from \eqref{eq:fP0} and \eqref{eq:fP}.  The group $\G_5$ is the stabilizer of the first contact point for $\I^3_1=\{P_0,\lsp{w}P_0,\lsp{\tau w} P_0\}$.  The group $\G_6$ is the stabilizer of the first contact point for $\I^2_1=\{P_0,\lsp{w}P_0\}$. The group $\G_7$ is the stabilizer of point for $D(\I^8)$, where \[\I^8=\{P_0,\conj{w},\conj{w \tau \sigma w},\conj{\tau^{-1} \sigma w},\conj{\tau \check{\sigma} w},\conj{\tau w \tau \sigma w},\conj{\tau^2 \check{\sigma} \sigma w},\conj{\epsilon w \xi^4 w}\}.\]
The group $\G_8$ is the stabilizer of the first contact point for $\I^3_2=\{P_0,\lsp{w}P_0,\conj{\sigma w}\}$.  The group $\G_9$ is the stabilizer of the first contact point for $\I^2_2=\{P_0,\conj{\xi}\}$.  The computation of first contact points in \cite[Proposition 6.5]{Yaspicard} concludes the proof.
\end{proof}
\end{subsection}
\end{section}
\bibliography{../../references}    
\end{document}